\documentstyle[11pt, letterpaper]{amsart}
\def\R{{\Bbb R}}

\def\Rn{{{\Bbb R}^n}}

\def\Rdvan{{{\Bbb R}^{2n}}}

\def\ro{{\rho}}
\def\la{{\lambda}}

\def\dz{{\zeta}}
\def\et{{\eta}}
\def\ta{{\tau}}
\def\al{{\alpha}}
\def\be{{\beta}}
\def\ga{{\gamma}}
\def\de{{\delta}}

\def\dde{{\Delta}}
\def\xdot{{\dot{x}}}
\def\x{{x}}

\def\aln{{\al_0}}
\def\an{{a_0}}
\def\un{{u_0}}
\def\vn{{v_0}}
\def\Un{{U_0}}
\def\Uj{{U_1}}
\def\Ud{{U_2}}
\def\Ga{{\Gamma}}
\def\gab{{G_{ab}}}
\def\xj{{x_1}}
\def\xd{{x_2}}
\def\g{{g}}
\def\gg{{G}}
\def\h{{h}}

\def\f{{f}}
\def\ff{{F}}
\def\d{{d}}
\def\ss{{S}}
\def\qq{{Q}}
\def\kk{{K}}
\def\zlomek#1#2{{\textstyle\frac{#1}{#2}}}

\begin{document}
\sloppy
\pagestyle{plain}
\setcounter{page}{1}
\frenchspacing
\title{The functional formulation of second-order ordinary differential equations}
\author{Petr Chl\'{a}dek}
\curraddr{Mathematical Institute of the Silesian University at Opava\\
\newline
\indent
Na Rybn\'\i\v{c}ku 1\\
746 01 Opava\\Czech Republic\\}
\makeatletter
\email{Petr.Chladek@math.slu.cz}
\makeatother
\begin{abstract}
\noindent
In this paper, the necessary and sufficient conditions in order that a smooth mapping $\ff(\ta,\al,\be,a,b)$
be a dependence of a complete solution $\x(\ta)$ of some second-order ordinary differential equation on Neumann
conditions $\x(\al)=a$, $\x(\be)=b$, $\al \neq \be$ are deduced. These necessary and sufficient conditions consist
of functional equations for $\ff$ and of a smooth extensibility condition. Illustrative examples are presented to
demonstrate this result. In these examples, the mentioned functional equations for $\ff$ are related to the
functional equations for geodesics, to Jensen's equation, to the functional equations for conic sections and
to Neuman's result for linear ordinary differential equations. 
\end{abstract}
\maketitle
\vspace{-1cm}
\section{Introduction}
\noindent
The relationship between functional and ordinary differential equations was studied by many authors, e.g., M. A.
Abdelkader \cite{abdelkader}, E. Castillo and R. Ruiz-Cobo \cite{castillo}, L. N. Eshukov \cite{eshukov}. We will
discuss here a functional formulation of second-order ordinary differential equations.
\par
In Section 2, we present the functional equations (\ref{8}), (\ref{9}) for a continuous function of many variables.
A solution for 1-dimensional linear case
is found. This solution is related to Neuman's result for linear second-order ordinary differential equations
\cite{neuman}.
\par
In Section 3, we investigate the relationship between the introduced functional equations and ordinary differential
equations in general. A notion of so-called N-solvable second-order ordinary differential equation is presented. We
prove that under given condition (\ref{10}) the functional equations are equivalent to the family of all smooth
N-solvable second-order ordinary differential equations.
\par
Section 4 contains examples. The introduced functional equations are related to Jensen's functional equation, 
to the functional equations for geodesics, and also to the functional equations for conic sections. 
\par
For the sake of simplicity, in this paper we consider smooth ordinary differential equations globally defined on
the Cartesian product and globally N-solvable.\\
\par
\section{Functional equations}
\noindent
Let $n$ be a natural number, $D=\{(\al,\be) \in \R^2 \mid \al=\be \}$ the diagonal of $\R^2$.
Let 
$\ff \colon \R \times (\R^2 \setminus D) \times \Rdvan \rightarrow \Rn$
be a continuous mapping.
We consider functional equations
\begin{equation}
\label{8}
\ff(\ta,\al,\be,a,b)=\ff(\ta,\ga,\de,\ff(\ga,\al,\be,a,b),\ff(\de,\al,\be,a,b)),
\end{equation}
\begin{equation}
\label{9}
\ff(\al,\al,\be,a,b)=a, \quad \ff(\be,\al,\be,a,b)=b.
\end{equation}
Let us start with an example.\\
\par
\noindent
{\bf Example 1} (linear case). Let $n=1$ and let $\ff$ be linear in the last two arguments.
Then $\ff$ is a solution of the functional equations (\ref{8}), (\ref{9}) if and only if there
exist functions $\xj,\xd\colon \R \rightarrow \R$ such that whenever $\al \neq \be$ the following statements hold
\begin{equation}
\label{l4}
\xj(\al)\xd(\be)-\xj(\be)\xd(\al) \neq 0,
\end{equation}
\begin{equation}
\begin{array}{c}
\label{l3}
\ff(\ta,\al,\be,a,b)=\qquad\qquad\qquad\\
\\
{\displaystyle
\frac{(\xj(\ta)\xd(\be)-\xj(\be)\xd(\ta))a+(\xj(\al)\xd(\ta)-\xj(\ta)\xd(\al))b} 
{\xj(\al)\xd(\be)-\xj(\be)\xd(\al)} } .
\end{array}
\end{equation}
\\
Let us check it.
Conditions (\ref{8}), (\ref{9}) can be obtained from (\ref{l4}), (\ref{l3}) by direct computation.
Conversely, put 
$\xj(\ta)=\ff(\ta,0,1,1,0)$,
$\xd(\ta)=\ff(\ta,0,1,0,1)$.
Then
\begin{equation}
\label{l5}
\mbox{det}\left( \begin{array}{ccc} 
\xj(\al) & \xd(\al) & a\\
\xj(\be) & \xd(\be) & b\\
\xj(\ta) & \xd(\ta) & \ff(\ta,\al,\be,a,b)
\end{array} \right)=0.
\end{equation}
Indeed, from linearity and from (\ref{8}) the third row is a linear combination of the first row and the second row
with coefficients $\ff(\ta,\al,\be,1,0)$ and $\ff(\ta,\al,\be,0,1)$.
Since from linearity and from (\ref{8}), (\ref{9})
$$
\left( \begin{array}{cc} 
\ff(0,\al,\be,1,0) & \ff(0,\al,\be,0,1) \\
\ff(1,\al,\be,1,0) & \ff(1,\al,\be,0,1) 
\end{array} \right)
\left( \begin{array}{cc} 
\xj(\al) & \xd(\al) \\
\xj(\be) & \xd(\be) 
\end{array} \right)=
\left( \begin{array}{cc} 
1 & 0 \\
0 & 1 
\end{array} \right) 
$$
holds,
we have the inequality (\ref{l4}). Finally, from (\ref{l5}) we can express $\ff$ in the form (\ref{l3}).\\
\par
\noindent
{\it Remark $1.$} Consider a second-order linear differential equation with linearly
independent solutions
$\xj,\xd$. For this case the condition (\ref{l5}), where $\ff$ is the solution satisfying (\ref{9}),
was presented (in a more general form) by F. Neuman
\cite{neuman} in 2003. In the following chapter we will investigate the relationship between
the functional equations (\ref{8}), (\ref{9}) and differential equations in general.\\
\par
\section{The relationship between functional and differential equations}
\noindent
Let $\f \colon \R \times \Rn \times \Rn \rightarrow \R^n$ be a mapping.
Consider a second-order equation
\begin{equation}
\label{1}
\ddot{\x}=\f(\ta,\x,\xdot)
\end{equation}
and Neumann conditions 
\begin{equation}
\label{4}
\x(\al)=a, \quad
\x(\be)=b,
\end{equation}
where $\al,\be \in \R$, $\al \neq \be$, $a,b \in \Rn$.
\par
The mapping $\ff$ will be called a {\it dependence of the solution on Neumann conditions} if and only if the
following equality holds
\begin{equation}
\label{7}
\ff(\ta,\al,\be,\x(\al),\x(\be))=\x(\ta)
\end{equation}
for any $\al \neq \be$, $\ta \in \R$ and for any complete solution 
$\x$ of the equation (\ref{1}). Notice that in this case $\mbox{Dom}(\x)=\R$.
We say that the dependence of the solution on Neumann conditions is {\it smooth} if and only if
$\ff$ is a smooth mapping.
The equation (\ref{1}) is called {\it N-solvable} if and only if there exists a unique complete solution
for any Neumann conditions (\ref{4}).
Consider the equation (\ref{1})
with integral conditions
\begin{equation}
\label{2}
\x(\al)=a, \quad
\int \limits _0^1 \xdot((1-\ga)\al+\ga\be)\,d\ga=v,
\end{equation}
where $\al,\be,\ga \in \R$, $a,v \in \R^n$.
The equation (\ref{1}) is called {\it I-solvable} if and only if there exists a unique complete solution for any
integral conditions.\\
\par
\noindent
{\bf Lemma 1.} The smooth equation (\ref{1}) is N-solvable if and only if it is I-solvable.\\
\par
\noindent
{\it Proof.} Let (\ref{1}) be N-solvable and smooth. For $\al \neq \be$ we can rewrite (\ref{2}) as
\begin{equation}
\label{5}
\x(\al)=a, \quad
\x(\be)=a+v(\be-\al).
\end{equation}
For $\al=\be$ we have (\ref{2}) in the form
\begin{equation}
\label{6}
\x(\al)=a, \quad
\xdot(\al)=v.
\end{equation}
The existence and uniqueness of the solution of (\ref{1}), (\ref{5}) (respectively, (\ref{1}),
(\ref{6})) follows from the assumption of N-solvability (respectively, from the uniqueness of the solution
of the Cauchy problem, see, e.g., \cite{hartman}).  Therefore (\ref{1}) is I-solvable.
\par
Conversely, let (\ref{1}) be I-solvable.
Then (\ref{1}) is N-solvable. This follows from (\ref{5}) immediately. $\Box$\\
\par
\noindent
{\bf Lemma 2.} The dependence $(\ta,\al,\be,a,v) \mapsto \x(\ta)$ of the solution $\x$ of the smooth
I-solvable equation (\ref{1}) on conditions (\ref{2}) is smooth on an open set containing $\R \times D
\times
\Rdvan$, where $D$ is the diagonal of $\R^2$.\\
\par
\noindent
{\it Proof.} Let $\ta \mapsto \x(\ta)$ be the solution of (\ref{1}) satisfying Cauchy conditions
$\x(\al)=a$, $\xdot(\al)=u$.
The dependence of $\x$ on $\ta,\al,a,u$ is smooth (see, e.g., \cite[chapter V]{hartman}). 
The mappings
$$
\h \colon (\al,\be,a,v,u) \mapsto \int \limits _0^1 \xdot((1-\ga)\al+\ga\be)\,d\ga-v,
$$
$$
\d \colon (\al,\be,a,v,u) \mapsto \mbox{det}
\left( \frac{\partial \h(\al,\be,a,v,u) }{\partial {u}} \right)
$$
are smooth too.
It remains to check that the dependence of $u$ on $\al,\be,a,v$ is smooth.
If we put $\be=\al$ in (\ref{2}), using Cauchy conditions we get $u=v$.
Since
$\h(\aln,\aln,\an,\vn,\vn)=0$ and $\d(\aln,\aln,\an,\vn,\vn)=1$ for
$\aln \in \R$, $\an,\vn \in \Rn$, we can apply \cite[Theorem 1.3.5. and Corollary 1.3.9.]{narasimhan} 
to the mapping $h$.
From this there exists a neighbourhood 
$\Uj \subset \R^{2n+2} \mbox{ of } (\aln,\aln,\an,\vn)$ and a neighbourhood $\Ud \subset \Rn \mbox{ of }
\vn$ such that for any
$(\al,\be,a,v)
\in
\Uj$ there is a unique $u=u(\al,\be,a,v) \in \Ud$, which satisfies 
$\h(\al,\be,a,v,u)=0$ and $(\al,\be,a,v) \mapsto u(\al,\be,a,v)$ is a smooth mapping.
Finally, the dependence of $\x$ on $\ta$ and on parameters $(\al,\be,a,v)$ from $\Uj$ is smooth. $\Box$\\ 
\par
\noindent
{\bf Theorem.} A mapping $\ff$ is the smooth dependence of the solution of some smooth
N-solvable equation (\ref{1})
on Neumann conditions if and only if the conditions (\ref{8}), (\ref{9}) hold
and a mapping
\begin{equation}
\begin{array}{c}
\label{10}
\R \times (\R^2 \setminus D) \times \Rdvan \ni (\ta,\al,\be,a,v) \mapsto \qquad\qquad\qquad \\
\qquad\qquad\qquad\qquad  \ff(\ta,\al,\be,a,a+v(\be-\al)) \in \Rn
\end{array}
\end{equation}
has a smooth extension $S$ to $\R^{2n+3}$.\\
\par
\noindent
{\it Proof of Theorem.} Assume that $\ff$ is the smooth dependence of the solution of some smooth 
N-solvable equation (\ref{1}) on Neumann conditions (\ref{4}). 
From (\ref{4}), (\ref{7}) we get 
$\ff(\ta,\al,\be,a,b)=\x(\ta)=\ff(\ta,\ga,\de,\x(\ga),\x(\de))$
for $\al \neq \be$, $\ga \neq \de$, which implies (\ref{8}).
If we put $\al=\ta$ (respectively, $\be=\ta$) in (\ref{7}), using (\ref{4}) we obtain first
equality in (\ref{9}) (respectively, second equality in (\ref{9})). 
We define a mapping
$\ss \colon (\ta,\al,\be,a,v) \mapsto \x(\ta)$,
where $\x$ is the solution of the smooth N-solvable equation (\ref{1}) satisfying conditions (\ref{2}). 
The mapping $\ss$ exists for any $(\ta,\al,\be,a,v) \in \R^{2n+3}$, this follows from Lemma 1.
Moreover, from (\ref{7}), (\ref{5}) we have 
$\ss(\ta,\al,\be,a,v)=\ff(\ta,\al,\be,a,a+v(\be-\al))$
for $\al \neq \be$.
Therefore $\ss$ is the extension of (\ref{10}) to $\R^{2n+3}$. The smoothness of $\ss$
follows from the smoothness of $\ff$ and from Lemma 2.
\par
Now we assume that (\ref{8}), (\ref{9}), (\ref{10}) hold.
With help of (\ref{10}) we calculate the limit $\de \rightarrow \ga$ in (\ref{8}) as
$$
\ss(\ta,\al,\be,a,v)=\ss\left(\ta,\ga,\ga,\ss(\ga,\al,\be,a,v),
\frac{\partial \ss(\ga,\al,\be,a,v) }{\partial {\ga}}  \right).
$$
Differentiating this with respect to $\ta$ two times and putting $\be=\al$, $\ga=\ta$
we see that $\ta \mapsto \ss(\ta,\al,\al,a,v)$ is the solution of (\ref{1}),
where
\begin{equation}
\label{12}
\left.
\f(\ta,\x,\xdot)=
\frac{\partial^{2} \ss(\ta,\ga,\ga,\x,\xdot) }{\partial {\ta^2}} \right |_{\ga=\ta}.
\end{equation}
From (\ref{9}), (\ref{10}) we obtain
$\ss(\al,\al,\al,a,v)=a$
and
$$
\left. \frac{\partial \ss(\ta,\al,\al,a,v)}{\partial \ta} \right |_{\ta=\al}=
\left. \left. \frac{\partial \ss(\ta,\al,\be,a,v)}{\partial \ta} \right |_{\ta=\be}
\right |_{\be=\al}=
$$
$$
=\left. \left(  \frac{\partial}{\partial \be} \left (\ss(\ta,\al,\be,a,v)
\left. \right |_{\ta=\be} \right )-
\left. \frac{\partial \ss(\ta,\al,\be,a,v)}{\partial \be} \right |_{\ta=\be}
\right) \right |_{\be=\al}=
$$
$$
=\left. \frac{\partial}{\partial \be} \left( \left. \ss(\ta,\al,\be,a,v)
\right |_{\ta=\be}-
\left. \ss(\ta,\al,\be,a,v) \right |_{\ta=\al} \right)
\right |_{\be=\al}=
$$
$$
=\left. \frac{\partial}{\partial \be} \left( a+v(\be-\al)-a \right) 
\right |_{\be=\al}=v.
$$
From (\ref{12}) $\f$ is smooth.
Let $\x$ be a complete solution of (\ref{1}), (\ref{12}).
From above we can see that the solutions $\ta \mapsto \ss(\ta,\al,\al,\x(\al),\dot{\x}(\al)) \mbox{ and }
\ta \mapsto \x(\ta)$ 
satisfy the same Cauchy condition. Therefore 
\begin{equation}
\label{14}
\ss(\ta,\al,\al,\x(\al),\dot{\x}(\al))=\x(\ta) 
\end{equation}
for any $\ta \in \R$
(see, e.g., \cite[chapter V]{hartman}).
Finally, as $\be \rightarrow \al$ in (\ref{8}), using (\ref{4}), (\ref{10}), (\ref{14})
we obtain $\ff(\ta, \ga, \de, \x(\ga), \x(\de))=\x(\ta)$, which is equivalent to (\ref{7}).
Hence $\ff$ is the dependence of the solution on Neumann conditions.
The smoothness of $\ff$ follows from (\ref{10}).
The uniqueness of the complete solution of (\ref{1}), (\ref{4}), (\ref{12}) follows from (\ref{7}).
Thus the equation (\ref{1}), (\ref{12}) is N-solvable. $\Box$\\
\par
\noindent
{\it Remark $2.$} From Theorem we can see that under the condition (\ref{10}) the functional equations (\ref{8}),
(\ref{9}) are equivalent to the family of all smooth N-solvable ordinary differential equations (\ref{1}).
It follows from the proof of Lemma 2 that every smooth differential equation (\ref{1}) is locally
I-solvable, i.e., for any $\aln \in \R$, $\an, \un \in \Rn$ there exists a neighbourhood 
$\Un \subset \R^{2n+1} \mbox{ of } (\aln,\an,\un)$ and a neighbourhood $\Uj \subset \R^{2n+2} \mbox{ of }
(\aln,\aln,\an,\un)$ such that the equation (\ref{1}) restricted on $\Un$
has a unique complete solution for any $(\al,\be,a,v) \in \Uj$. 
Moreover, from the proof of Lemma 1, every smooth equation (\ref{1}) is locally N-solvable.
Thus the functional equations (\ref{8}), (\ref{9}) hold locally for any smooth equation (\ref{1}), 
but $\mbox{Dom}(F)$ is more complicated in such a case.\\ 
\par
\section{Examples}
\par 
\noindent
{\bf Example 2} (free fall). Let us consider the free fall equation
$\ddot{\x}=\g$.
It is easy to see, that for this equation we have
\begin{equation}
\begin{array}{c}
\label{f1}
\ff(\ta,\al,\be,a,b)={\displaystyle \frac{(\ta-\be)a+(\al-\ta)b}{\al-\be}}+\zlomek{1}{2}\g(\ta-\al)(\ta-\be)
\end{array}
\end{equation}
and a smooth extension (\ref{10}) has the form
$$
(\ta,\al,\be,a,v) \mapsto a+v(\ta-\al)+\zlomek{1}{2}\g(\ta-\al)(\ta-\be).
$$
\medskip
\par
\noindent
{\bf Example 3} (geodesics). A geodesic can be viewed as a solution of a second-order differential equation.
Thus this is the special case of the general situation presented above. Let $\Ga$ be a smooth linear connection on a
manifold
$M$. Let $U \subset M$ be an open convex
set and let $a,b \in U$. By $\gab$ we denote the geodesic of $\Ga$ satisfying $\gab(0)=a$, $\gab(1)=b, \gab([0,1])
\subset U$. The existence of $U$ and uniqueness of $\gab$ is guaranteed by the Whitehead local convex lemma 
(see, e.g., \cite{kobayashi} or \cite{whitehead}). 
We define a mapping
$\gg$ in the following way
$$\gg \colon (a,b,\ro) \mapsto \gab(\ro).$$
Since there exist constant geodesics, we have
\begin{equation}
\label{g0}
\gg(a,a,\ro)=a.
\end{equation}
Let $\x$ be a geodesic satisfying 
$\x(\al)=a \mbox{, } \x(\be)=b$. 
Since every affine reparametrization of a geodesic is a geodesic again, we have
$\x(\ta)=\gg(a,b,\frac{\ta-\al}{\be-\al})$.
If there exists a dependence of the solution on Neumann conditions, then
from the previous equality and from (\ref{7}) we obtain
\begin{equation}
\label{g-1}
\ff(\ta,\al,\be,a,b)=\gg\left(a,b,\frac{\ta-\al}{\be-\al}\right),
\end{equation}
where $\mbox{Dom}(\gg)=\Rn \times \Rn \times \R$ and $\mbox{Codom}(\gg)=\Rn$.
The following functional equations for geodesics were deduced by L. Klapka \cite{klapka} in 2000
\begin{equation}
\label{g2}
\gg(a,b,0)=a, \quad \gg(a,b,1)=b,
\end{equation}
\begin{equation}
\label{g3}
\gg(a,b,(1-\ro)\dz+\ro \et)=\gg(\gg(a,b,\dz),\gg(a,b,\et),\ro).
\end{equation}
Under the conditions (\ref{g0}), (\ref{g-1}),
the statements (\ref{g2}), (\ref{g3}) follow from (\ref{8}), (\ref{9})
for any $\dz,\et,\ro \in \R$. \\
\par
\noindent
{\bf Example 4} (Jensen's functional equation). Consider a mapping $\gg$ from Example 3. Let $\gg$ be linear in
the first two arguments. By putting $\dz=1$, $\et=0$ in (\ref{g2}), (\ref{g3}) we obtain
$\gg(a,b,1-\ro)=\gg(b,a,\ro)$.
For any $\ro \in \R$ we define a mapping $\qq(\ro)$ in the following way
$$\qq(\ro) \colon \Rn \ni b \mapsto \gg(0,b,\ro) \in \Rn.$$
From this and from linearity we get
$\gg(a,b,\ro)=\qq(1-\ro)(a)+\qq(\ro)(b)$
and using (\ref{g0}) we obtain
$\qq\left(\zlomek{1}{2}\right)(a)=\zlomek{1}{2}a$.
Finally, by substituting $a=0$, $\ro=\frac{1}{2}$ into (\ref{g3}) and using the previous equalities we get
Jensen's functional equation
$$
\qq\left(\frac{\dz+\et}{2}\right)=\frac{\qq(\dz)+\qq(\et)}{2}.
$$\\
\par
\noindent
{\bf Example 5} (Conic sections). Let $n=1$. Let $\ff$ be continuous,
linear inhomogeneous in the last two arguments, autonomous, i.e.,
\begin{equation}
\label{c7}
\ff(\ta,\al,\be,a,b)=\ff(\ta+\la,\al+\la,\be+\la,a,b)
\end{equation}
for any $\la \in \R$ and nondissipative, i.e.,
$\ff(\ta,\al,\be,a,b)=\ff(-\ta,-\al,-\be,a,b)$.
Under these assumptions we can find a solution of (\ref{8}), (\ref{9}).
Let $\al,\be,a,b$ be fixed. We denote $\x(\ta)=\ff(\ta,\al,\be,a,b)$.
From (\ref{8}), (\ref{c7}) we get
$$
\x(\ta)=\ff(0,\dde,2\dde,\x(\ta+\dde),\x(\ta+2\dde)),
$$
$$
\x(\ta+3\dde)=\ff(0,-\dde,-2\dde,\x(\ta+2\dde),\x(\ta+\dde))
$$
for any $\ta \in \R$, $\dde \in \R \setminus \{0\}$. By using the inhomogeneous linearity and the nondissipativity
we see that the difference
$\x(\ta+3\dde)-\x(\ta)$ depends on $\x(\ta+2\dde)-\x(\ta+\dde)$ linearly and a coefficient depends on $\dde$ only.
Similarly, $\x(\ta+4\dde)-\x(\ta+\dde)$ depends on $\x(\ta+3\dde)-\x(\ta+2\dde)$ linearly and the coefficient
is the same.
Therefore
$$
(\x(\ta+4\dde)-\x(\ta+\dde))(\x(\ta+2\dde)-\x(\ta+\dde))= \qquad\qquad\qquad
$$
$$
\qquad\qquad\qquad\qquad\qquad(\x(\ta+3\dde)-\x(\ta))(\x(\ta+3\dde)-\x(\ta+2\dde)).
$$
This is the functional equation for characterization of conic sections presented by A. Angelesco
(see \cite[section 2.5.2.]{aczel} or \cite{angelesco}) in 1922.
All its continuous solutions are known and
these solutions are real analytic. Differentiating this equation with respect to $\dde$ five times and using the
limit
$\dde \rightarrow 0$ we get an equation
$$
\xdot\,\ddddot{\x}=\dddot{x}\,\ddot{\x},
$$
which is equivalent to 
$$
\frac{\mbox{d}}{\mbox{d}\ta}\left(\frac{\dddot{\x}(\ta)}{\xdot(\ta)}\right)=0.
$$
From this
$\dddot{x}=\kk\xdot$
and by integrating we obtain 
$\ddot{\x}=\kk\x+g$,
where $K,g \in \R$. According to N-solvability $K \geq 0$.
Putting $\kk=k^2$ we get
$$
f(\ta,\x,\xdot)=k^2\x+g.
$$
From (\ref{1}), (\ref{7}) 
$$
\ff(\ta,\al,\be,a,b)=
$$
$$
\frac{(k^2a+g)\mbox{sinh}(k(\ta-\be))+(k^2b+g)\mbox{sinh}(k(\al-\ta))-g\,\mbox{sinh}(k(\al-\be))}
{k^2\, \mbox{sinh}(k(\al-\be))},
$$
where $k,g \in \R$, $k \neq 0$. For $k=0$ we get the free fall (\ref{f1}).\\
\par
\noindent
{\bf Acknowledgment.} The author would like to thank L. Klapka for his great help during the work.
This research was supported by the grant MSM 192400002 from Czech Ministry of Education.
\par
\renewcommand{\refname}{References}

\par
\end{document}